\documentclass[11pt]{amsart}
\usepackage{amsfonts,amsmath,amssymb,enumerate}

\newcommand{\bpf}[1][Proof]{{\noindent {\sc #1: }}}
\newcommand{\epf}{{{\hfill $\Box$ \smallskip}}}
\newtheorem{lemma}{Lemma}
\newtheorem{theorem}{Theorem}

\newtheorem{remark}{Remark}

\newcommand{\Fc}{\mathcal{F}}
\newcommand{\Pp}{\mathsf{P}}
\newcommand{\R}{\mathbb{R}}

\newcommand{\ONE}{{\bf 1}}

\newcommand{\eps}{\varepsilon}

\newcommand{\E}{\mathsf{E}}
\newcommand{\Law}{\mathop{\mathsf{Law}}}

\begin{document}
\title{On Gumbel limit for the length of reactive paths}
\author{Yuri Bakhtin}
\address{686 Cherry Street, School of Mathematics, Georgia Institute of Technology, Atlanta, GA, 30332-0160}
\email{yuri.bakhtin@gmail.com}
\date{}
\thanks{The author was partially supported by the NSF CAREER Award DMS-0742424.}
\begin{abstract} We give a new proof of the vanishing noise limit theorem for exit times of 1-dimensional
diffusions conditioned on exiting through a point separated from the starting point by a potential wall. 
We also prove a scaling limit for exit location in a model 2-dimensional situation.
\end{abstract}


\maketitle
\section{Introduction}
A random variable $Z$ has standard Gumbel distribution if 
\begin{equation*}
\Pp\{Z\le x\}=\Lambda(x)=e^{-e^{-x}},\quad x\in\R.
\end{equation*}
The Gumbel distribution is mostly well-known as a limit law in the extreme value theory, see, e.g.,
\cite[Chapter 1]{deHaan:MR2234156}. 

Surprisingly, it also
appears as a limiting distribution (in the limit of vanishing noise)
for normalized
exit times of diffusions conditioned on unlikely exit locations. The first result of
this kind was obtained in~\cite{Day:MR1110156}, where the phenomenon of cycling for exits through a characteristic repelling boundary was
discovered and explained via the asymptotics of exit times, see~\cite{BG_periodic2}
for most recent progress on cycling and a bibliography.
 
 In a recent paper~\cite{ALEA} (see also references therein for related results) a similar
result for conditional exit times was obtained for the case where the prescribed exit location for the diffusion is separated
from the starting point by a potential wall.

Namely, consider a strong solution $X_\eps$ of the following SDE driven by additive white noise of small intensity
$\eps>0$ and a smooth vector field $b$ on an interval $[q_-,q_+]$ containing $0$:
\begin{align}
\label{eq:basic-sde}
dX_\eps(t)&=b(X_\eps(t))dt+\eps dW(t),\\
\label{eq:initial-condition}
X_\eps(0)&=x_0,
\end{align}
Here $x_0\in(q_-,0)$, the drift $b$ satisfies $b(0)=0$, $b'(0)>0$,
$b(x)<0$ for $x\in[q_-,0)$, and $b(x)>0$ for $x\in(0,q_+]$, and $W$ is a standard Wiener process on a complete probability
space $(\Omega,\Fc,\Pp)$.

We are interested in the conditional distribution of the first exit time \[\tau_\eps=\inf\bigl\{t\ge 0: X_\eps(t)\in\{q_-,q_+\}\bigr\},\]
conditioned on the event $C_\eps(q_-,x_0,q_+)=\{X_\eps(\tau_\eps)=q_+\}$. This is an extremely improbable event for small $\eps$ 
since for $C_\eps(q_-,x_0,q_+)$ to happen, the process $X_\eps$ has to travel against the drift. The resulting diffusion trajectories are often called reactive paths.

\begin{theorem}[\cite{ALEA}]\label{th:ALEA} If, in addition to the conditions given above, $q_+>0$, then there are constants $c_1,c_2,c_3$ such that
as $\eps\to0$, 
\[
\Law\left[\tau_\eps-c_1\ln\frac{1}{\eps}\ \Bigr|\ C_\eps(q_-,x_0,q_+)\right]\ \Rightarrow\   \Law\left[c_2Z+c_3\right],
\]
where $Z$ is a standard	Gumbel random variable and ``$\Rightarrow$'' denotes weak convergence
of probability measures.
\end{theorem}

The core of the argument in~\cite{ALEA} was based on asymptotic analysis of solutions of second-order differential
equations describing the moment generating functions of exit times. This is a powerful general method, and 
in~\cite{ALEA} it was also used to study other interesting settings. Unfortunately,
 it does not seem to answer why it is natural for the max-stable Gumbel distribution to appear in this context.
The path-based arguments for slightly different situations in~\cite{Day:MR1110156} and~\cite{Day:MR1376805} are more elucidating,
and it would be natural to supplement the argument of~\cite{ALEA} by a more probabilistic one.

In~\cite{Bakhtin-Gumbel}
we observed that the connection between the extreme value theory and exit times is provided by the theory of
residual life times, see~\cite{residual:MR0359049}. However, to make our point in that paper we chose to work with an initial
condition $x_0$ depending on $\eps$, namely, we set $x_0=-\eps a$ and considered two consecutive limit transitions
instead of one: first $\eps\to 0$, and then $a\to+\infty$. One of the goals of the present paper is to show that purely
probabilistic reasoning 
based on basic properties of the Wiener process is, in fact, not much harder for the original problem with fixed $x_0\in(q_-,0)$
and only one limiting transition $\eps\to 0$.

As in~\cite{Bakhtin-Gumbel}, we make a simplifying assumption $b(x)=\lambda x$ for a constant $\lambda>0$ 
(in general, one has to approximate the nonlinear case with a linearization). We will prove the following result with
purely probabilistic tools:
\begin{theorem}[\cite{ALEA}] 
\label{th:main} Let $Z$ denote a standard Gumbel random variable.
If $q_+>0$, then
\begin{equation}
\label{eq:main-theorem-exit-right}
\Law\left[\tau_\eps-\frac{2}{\lambda}\ln\frac{1}{\eps}\Bigr|\ C_\eps(q_-,x_0,q_+)\right]\Rightarrow\Law\left[\frac{\ln (2\lambda q_+|x_0|)}{\lambda}+\frac{1}{\lambda}Z\right],\quad \eps\to 0.
\end{equation}
\end{theorem}
We stress that this result per se is not new. In~\cite{ALEA} it is given as Proposition~3.1  and plays the central role in the proof of
Theorem~\ref{th:ALEA} of which it is a specific case. However, our point is to give a proof based on properties of the Wiener process and to stress the connection to residual life times.
We give our proof in Section~\ref{sec:proof-of-Gumbel}. Our method is somewhat similar to the proof of the following statement which is a version
of a result from~\cite{Day:MR1110156} (where, in fact, a more general nonlinear version is considered) for the case of characteristic boundary, i.e. $q_+=0$:  
\begin{theorem}[\cite{Day:MR1110156}]
\label{th:Day}
\begin{equation}
\label{eq:main-theorem-exit-0}
\Law\left[\tau_{\eps}-\frac{1}{\lambda}\ln\frac{1}{\eps}\Bigr|\ C_\eps(q_-,x_0,0)\right]\Rightarrow\Law\left[\frac{\ln(x_0^2\lambda)}{2\lambda}+\frac{1}{2\lambda}Z\right],
\quad \eps\to 0.
\end{equation}
\end{theorem}
For completeness we give a proof of this result in Section~\ref{sec:proof-of-Gumbel}, too.

We notice that Theorem~\ref{th:main} fits nicely with the results of~\cite{Day:MR1376805,Bakhtin:MR2411523,Bakhtin:MR2800902} 
which for the linear 1-dimensional situation with $q_+>0$ and $x_0=0$ take the following form:
\begin{theorem}[\cite{Day:MR1376805}] \label{th:exit-starting-at-zero}
Let $q_+>0$. Then
\begin{equation}
\label{eq:exit-starting-at-zero}
\Law\left[\tau_{\eps}-\frac{1}{\lambda}\ln\frac{1}{\eps}\Bigr|\ C_\eps(q_-,0,q_+)\right]\Rightarrow\Law\left[\frac{\ln q_+}{\lambda}+\frac{\ln(2\lambda)}{2\lambda}+\frac{1}{\lambda}\Theta\right],
\end{equation}
where $\Theta=-\ln|N|$ and $N$ is a standard Gaussian random variable.
\end{theorem}
The connection is the following: to realize the exit described in Theorem~\ref{th:main} the solution must first 
reach~0 and then reach~$q_+$ from~0. The strong Markov property allows us to conclude that the limit
in~\eqref{eq:main-theorem-exit-right} is equal to the convolution of limits in
~\eqref{eq:main-theorem-exit-0} and~\eqref{eq:exit-starting-at-zero}. This can be checked independently by a
straightforward computation of the Laplace
transforms or characteristic functions of the right-hand sides of~\eqref{eq:main-theorem-exit-right},~\eqref{eq:main-theorem-exit-0}, and ~\eqref{eq:exit-starting-at-zero}, and invoking the duplication formula for the gamma function.

\bigskip

The second goal of this paper is to consider systems with small noise in high dimensions and
use Theorem~\ref{th:main} to obtain a vanishing noise scaling limit 
for conditional  distributions of exit locations. Namely, we are interested in
exit from a neighborhood of a saddle point of the driving vector field under conditioning on a rare event that the solution
travels through the potential wall. 

In the unconditioned situation, scaling limits for exit distributions for neighborhoods of saddle points have been obtained in
~\cite{Bakhtin:MR2800902},~\cite{with-Almada:MR2802310}, and~\cite{with-Almada:MR2739004}.
These results are crucial for in the theory of heteroclinic networks under small noise in~\cite{Bakhtin:MR2800902} (see also~\cite{Bakhtin:MR2731621} for a less technical exposition). 

The only known to us results on scaling limits for exit distributions under conditioning on unlikely exit locations in multiple dimensions 
are those of~\cite{Bakhtin-Swiech} where a case with no invariant sets (such as critical points) in the domain under several
additional technical conditions is considered. 
The analysis in~\cite{Bakhtin-Swiech} is based on using Doob's $h$-transform to reduce the problem to
the Levinson case studied in~\cite{with-Almada:MR2739004}. However, the details are involved and require obtaining
a new gradient estimate for solutions of nonlinear elliptic equation.

Here we study only a model 2-dimensional situation with a linear saddle and additive isotropic noise that can be addressed using the
asymptotics provided by Theorem~\ref{th:main}, although we conjecture that the behavior is similar in the general nonlinear case of tunneling through a potential wall near a saddle point.

Let us consider the following stochastic dynamics in two dimensions:
\begin{align*}
dX^1_\eps(t)&=\lambda X^1_\eps(t)dt+\eps dW^1(t),\\
dX^2_\eps(t)&=-\mu X^2_\eps(t)dt+\eps dW^2(t),
\end{align*}
where $\mu$ and $\lambda$ are positive constants, and $W^1,W^2$ are two independent Wiener processes. The drift $b(x_1,x_2)=(\lambda x_1, -\mu x_2)
=-\nabla(-\lambda x_1^2+\mu x_2^2)/2$
of this model system is linear and the origin is a saddle critical point with first axis serving as the unstable manifold and second axis
as the stable one. 

Let us assume that the initial conditions $(X^1_\eps(0),X^2_\eps(0))$ for this system satisfy the following scaling limit: 
$X^1_\eps(0)=x_1$, where $x_1<0$ is a constant, and $X^2_\eps=\eps^\alpha\xi_\eps$ for some  $\alpha\in\R$ and a family of random
variables $(\xi_\eps)_{\eps>0}$ independent of the driving Wiener processes $W^1,W^2$ and  converging in distribution as $\eps\to0$ to a random variable $\xi_0$ that is not concentrated at~$0$.

Let us take two numbers $q_-,q_+$ satisfying $q_-<x_1<0<q_+$ and define the time $\tau_\eps$ as the first time the solution exits from
the strip $(q_-,q_+)\times \R$, i.e.,  $\tau_\eps=\inf\{t\ge 0:\ X^1_\eps(t)\in\{q_-,q_+\}\}$. We would like to study the exit distribution
of the process $X_\eps=(X^1_\eps,X^2_\eps)$ conditioned on $C_\eps=\{X^1_\eps(\tau_\eps)=q_+\}$, i.e., on the rare event that the solution makes it
over the potential wall to the other side of the saddle. 

\begin{theorem}\label{th:exit-scaling} Let $V,N,\xi$ be independent random variables such that $V$ is standard exponential, i.e., $\Pp\{V\ge x\}=e^{-x}$ for $x\ge 0$, 
$N$ is standard Gaussian, and $\xi$ is distributed as $\xi_0$. Let $\beta=1\wedge (2\mu/\lambda+\alpha)$. Then
\[
\Law\left[\frac{X^2_\eps(\tau_\eps)}{\eps^\beta}\Bigr|\ C_\eps\right]\Rightarrow\Law
\left[\ONE_{\beta=\frac{2\mu}{\lambda}+\alpha} \frac{V^{\mu/\lambda }\xi}{(2\lambda q_+|x_1|)^{\mu/\lambda}}
+\ONE_{\beta=1}\frac{N}{\sqrt{2\mu}}\right],
\quad \eps\to 0.
\]
\end{theorem}

In particular, we observe effects similar to those described in~\cite{Bakhtin:MR2800902} and~\cite{Bakhtin:MR2731621}. The character
of the limiting exit distribution under rescaling depends not only on the scaling exponent $\alpha$ of the distribution of the initial
condition, but also on the ratio of contraction  and expansion rates. 

In particular, if contraction is strong enough $(2\mu/\lambda+\alpha>1)$,
then at the exit the system asymptotically forgets the initial condition, and the exit distribution is asymptotically centered Gaussian and
scales as $\eps$, the magnitude of the noise. 

However, if the contraction is not strong enough $(2\mu/\lambda+\alpha<1)$, then the 
Gaussian term disappears and the dependence of the limiting exit distribution on the initial distribution is nontrivial. In particular,
this can lead to the following memory effect: if the distribution of $\xi$ is supported by the positive (or negative) semi-axis, then the 
limiting exit distribution will also retain this strong asymmetry.

 There is also an intermediate case $(2\mu/\lambda+\alpha=1)$ where both Gaussian and non-Gaussian terms contribute to the limiting distribution.

 \medskip
 
We give the proof of Theorem~\ref{th:exit-scaling} in Section~\ref{sec:scaling}.  It is clear from the proof that it is also easy to obtain a generalization of Theorem~\ref{th:exit-scaling} for higher dimensions.

\section{Proof of Theorems~\ref{th:main} and~\ref{th:Day}}\label{sec:proof-of-Gumbel}
\subsection{Asymptotic equivalence of families of events}
More than once in this section we will need to replace conditioning on one event by conditioning on another event. To that end it
is convenient to introduce a notion that will allow to control the change in the conditional probability.

We will say that families of events $(A_\eps)_{\eps>0}$ and $(B_\eps)_{\eps>0}$ are asymptotically equivalent as $\eps\to0$ if
$\Pp(A_\eps)>0$ for sufficiently small $\eps$ and
\begin{equation}
\label{eq:equivalence-of-events}
\lim_{\eps\to 0}\frac{\Pp(A_\eps\triangle B_\eps)}{\Pp(A_\eps)}=0.
\end{equation}
The roles of $A_\eps$ and $B_\eps$ are not symmetric in~\eqref{eq:equivalence-of-events}. However,
if~\eqref{eq:equivalence-of-events} holds, then $\Pp(B_\eps)=\Pp(A_\eps)(1+o(1))$, as $\eps\to 0$, and 
we see that $\Pp(B_\eps)>0$ for sufficiently small $\eps$ and
\begin{equation}
\label{eq:equivalence-of-events2}
\lim_{\eps\to 0}\frac{\Pp(A_\eps\triangle B_\eps)}{\Pp(B_\eps)}=0.
\end{equation}
Similarly,~\eqref{eq:equivalence-of-events} is implied by~\eqref{eq:equivalence-of-events2}, so the roles
of $A_\eps$ and $B_\eps$ are, in fact, interchangeable.
\begin{lemma}\label{lem:weak-convergence-under-equivalent-sequences} Let $(Y_\eps)_{\eps>0}$ be a family of random variables and assume that for a family of events
$(A_\eps)_{\eps>0}$, $\Law[Y_\eps|A_\eps]\Rightarrow \mu$  as
$\eps\to0$ for  a measure $\mu$.
If $(B_\eps)_{\eps>0}$ is a family of events asymptotically equivalent to $(A_\eps)_{\eps>0}$, then 
$\Law[Y_\eps|A_\eps]\Rightarrow \mu$  as
$\eps\to0$.
\end{lemma}
\bpf
We are given that for any continuous bounded function $f:\R\to\R$,
\begin{equation}
\label{eq:conditional-convergence}
\frac{\E[ f(Y_\eps)\ONE_{A_\eps}]}{\Pp(A_\eps)}\to \int_{\R}f(x)\mu(dx),\quad \eps\to 0.
\end{equation}
Therefore, for any continuous bounded function $f:\R\to\R$,
\begin{align*}
\frac{\E[ f(Y_\eps)\ONE_{B_\eps}]}{\Pp(B_\eps)}&=\frac{\E[ f(Y_\eps)\ONE_{A_\eps}]+\E[ f(Y_\eps)\ONE_{B_\eps\setminus A_\eps}]-
\E[ f(Y_\eps)\ONE_{A_\eps\setminus B_\eps}]}{\Pp(A_\eps)(1+o(1))}
\\&=
\frac{\E[ f(Y_\eps)\ONE_{A_\eps}]+o (\Pp(A_\eps))}{\Pp(A_\eps)(1+o(1))}\to \int_{\R}f(x)\mu(dx),\quad \eps\to 0.
\end{align*}
and the proof is completed.
\epf

\subsection{The representation of solution and reflection principle for an auxiliary process}
The strong solution of~\eqref{eq:basic-sde}--\eqref{eq:initial-condition} is given by variation of constants:
\begin{equation}
X_\eps(t)=e^{\lambda t}\left(x_0+\eps U(t)\right),\quad t\ge 0,
\label{eq:solution-Duhamel}
\end{equation}
where 
\[
U(t)=\int_0^te^{-\lambda s}dW(s),\quad t\ge 0.
\]
It will be useful for us to introduce
\[
D_\eps(z,r)=\left\{\eps U(t)= z\ \text{for some}\ t\in[0,r)\right\},\quad z\in\R,\ r\in[0,\infty],
\]
(in particular, $D_\eps(|x_0|,\infty)$ is the event that $X_\eps$ ever reaches $0$) 
and prove a representation for probabilities of these events. Let
\[
\Phi(x)=\frac{1}{\sqrt{2\pi}}\int_{-\infty}^x e^{-y^2/2}dy,\quad x\in\R,
\]
be the standard Gaussian distribution function.
\begin{lemma}\label{lem:calculating-prob-of-D} For any $z\in\R$, $r\in[0,\infty]$
\[
\Pp(D_\eps(z,r))=2\left(1-\Phi\left(\frac{|z|\sqrt{2\lambda}}{\eps\sqrt{1-e^{-2\lambda r}}}\right)\right).
\]
\end{lemma}
\bpf  Notice that the process $U$ admits a representation
\begin{equation*}
U(t)=B\left(\frac{1-e^{-2\lambda t}}{2\lambda}\right),\quad t\ge 0,
\end{equation*}
for a standard Wiener process $B$. 
Therefore,
\[
D_\eps(z,r)=\left\{\eps B(s)= z\ \text{for some}\ s\in\left[0,\frac{1-e^{-2\lambda r}}{2\lambda}\right)\right\},\quad r\in[0,\infty],
\]
so the the reflection principle for the Wiener process implies
\begin{align*}
\Pp(D_\eps(z,r))&=\Pp\left\{\sup_{s\in \left[0,\frac{1-e^{-2\lambda r}}{2\lambda}\right)}B(s)\ge\frac{|z|}{\eps}\right\}
=2\Pp\left\{B\left(\frac{1-e^{-2\lambda r}}{2\lambda}\right)\ge\frac{|z|}{\eps}\right\},
\end{align*}
and the lemma follows.\epf



\subsection{Proof of Theorem~\ref{th:Day}}
\label{sec:exit-to-0} This is essentially the same proof as in~\cite{Day:MR1110156} given here for completeness.
Let us define 
\[
\tau_{0,\eps}=\inf\bigl\{t\ge 0: X_\eps(t)=0\bigr\}=\inf\bigl\{t\ge 0: x_0+\eps U(t)=0\bigr\}\in[0,\infty].
\]
We have $\tau_{0,\eps}\ge \tau_{\eps}$ and $C_\eps:=C_\eps(q_-,x_0,0)=\{\tau_\eps=\tau_{0,\eps}\}$.
Let us denote 
\begin{equation}\label{eq:D_eps}
D_{\eps}=D_\eps(|x_0|,\infty)=\{\tau_{0,\eps}<\infty\}.
\end{equation}
Our goal is to replace conditioning on $C_\eps$ by conditioning on $D_\eps$.
\begin{lemma}\label{lem:C-similar-to-D}
Families $(C_\eps)_{\eps>0}$ and $(D_\eps)_{\eps>0}$ are asymptotically equivalent as $\eps\to 0$.
\end{lemma}
\bpf First we notice that $\Pp(C_\eps \setminus D_\eps)=0$. Next, using the strong Markov property, we obtain
\[
\frac{\Pp( D_\eps\setminus C_\eps)}{\Pp(D_\eps)}=
\frac{\Pp\{X_\eps(\tau_\eps)=q_-\}\Pp(D_\eps(|q_-|,\infty))}{\Pp(D_\eps(|x_0|,\infty))}.
\]
Now applying Lemma~\ref{lem:calculating-prob-of-D} and the standard estimate 
\begin{equation}
\label{eq:Gaussian-tail-estimate}
1-\Phi(x)\sim \frac{1}{x\sqrt{2\pi}}e^{-x^2/2},\quad x\to\infty,
\end{equation}
to $\Pp(D_\eps(|q_-|,\infty))$ and $\Pp(D_\eps(|x_0|,\infty))$,
we complete the proof. \epf

We turn now to  computing the conditional distribution 
of $\tau_{0,\eps}$ given $D_\eps$. Using Lemma~\ref{lem:calculating-prob-of-D} and Gaussian tail asymptotics~\eqref{eq:Gaussian-tail-estimate}, we obtain for any $r\in\R$:
\begin{align*}
&\Pp\left(\tau_{0,\eps}-\frac{1}{\lambda}\ln\frac{1}{\eps}< r\Bigr|\ D_\eps\right)=
\Pp\left(\tau_{0,\eps}<\frac{1}{\lambda}\ln\frac{1}{\eps}+ r\Bigr|\ D_\eps\right)
\\ =& \frac{\Pp\left(D_\eps\left(|x_0|,\frac{1}{\lambda}\ln\frac{1}{\eps}+r\right)\right)}{\Pp(D_\eps(|x_0|,\infty))}
= \frac{1-\Phi\left(\frac{|x_0|\sqrt{2\lambda}}{\eps\sqrt{1-e^{-2\lambda (\frac{1}{\lambda}\ln\frac{1}{\eps}+r)}}}\right)}{1-
\Phi\left(\frac{|x_0|\sqrt{2\lambda}}{\eps}\right)}
\\ =&\frac{1-\Phi\left(\frac{|x_0|\sqrt{2\lambda}}{\eps\sqrt{1-e^{-2\lambda r}\eps^2}}\right)
}{1-\Phi\left(\frac{|x_0|\sqrt{2\lambda}}{\eps}\right)}
\sim
\frac{\sqrt{1-e^{-2\lambda r}\eps^2}\cdot e^{-\frac{x_0^2\lambda}{\eps^2(1-e^{-2\lambda r}\eps^2)}}}{
e^{-\frac{x_0^2\lambda}{\eps^2}}}
\\
\to & e^{-x_0^2\lambda e^{-2\lambda r}},\quad \eps\to0.
\end{align*}
The right-hand side is the distribution function of $\frac{\ln(x_0^2\lambda)}{2\lambda}+\frac{1}{2\lambda}Z$. Lemma~\ref{lem:C-similar-to-D}
allows us to apply Lemma~\ref{lem:weak-convergence-under-equivalent-sequences} to $Y_\eps=\tau_{0,\eps}-\frac{1}{\lambda}\ln\frac{1}{\eps}$
and families $(D_\eps)_{\eps>0}$ and $(C_\eps)_{\eps>0}$. Since $\tau_\eps=\tau_{0,\eps}$ on $C_\eps$,~\eqref{eq:main-theorem-exit-0} follows.
\epf

\subsection{Proof of Theorem~\ref{th:main}} Let us recall that the solution $X_\eps$ is given by~\eqref{eq:solution-Duhamel}. We need to study the exit time $\tau_{\eps}$ conditioned on
the event $C_\eps(q_-,x_0,q_+)=\{X_\eps(\tau_\eps)=q_+\}$. From now on we use $C_\eps$ as a shorthand for
$C_\eps(q_-,x_0,q_+)$. In this proof we will change the conditioning twice.

Let us define
\[
\theta_\eps=\inf \left\{t\ge 0:\ X_\eps(t)=e^{\lambda t}\left(x_0+\eps U(t)\right)=q_+\right\}\in(0,\infty],
\]
and $E_\eps=\{\theta_\eps<\infty\}$. Notice that $C_\eps=\{\theta_\eps=\tau_\eps\}$.
\begin{lemma}\label{lem:C-similar-to-E} Families $(C_\eps)_{\eps> 0}$ and $(E_\eps)_{\eps> 0}$ are asymptotically equivalent as $\eps\to 0$.
\end{lemma}
\bpf The proof repeats the proof of Lemma~\ref{lem:C-similar-to-D}.
 \epf\medskip

 We will need one more asymptotic equivalence statement since it is most convenient to work with events $F_\eps=\{x_0+\eps U(\infty)>0\}.$
\begin{lemma} \label{lem:E-similar-to-F} Families $(E_\eps)_{\eps> 0}$ and $(F_\eps)_{\eps> 0}$ are asymptotically equivalent as $\eps\to 0$.
\end{lemma}
\bpf This proof is also very similar to that of Lemma~\ref{lem:C-similar-to-D}. First we notice that $\Pp(F_\eps\setminus E_\eps)=0$,
and then we use the strong Markov property to compute
\[
\frac{\Pp(E_\eps \setminus F_\eps)}{\Pp(E_\eps)}
=\Pp(F_\eps^c|E_\eps)=\Pp(D_\eps(q_+,\infty))\to 0,\quad \eps\to 0.
\]
\epf

We have $F_\eps\subset E_\eps$, i.e., on $F_\eps$ we have
\begin{equation}
e^{\lambda \theta_\eps}(x_0+\eps U(\theta_\eps))=q_+.
\label{eq:equation-for-theta}
\end{equation}
Therefore, on $F_\eps$ we have
\begin{align*}
\theta_\eps&=\frac{1}{\lambda}\ln\frac{q_+}{x_0+\eps U(\theta_\eps)}
         =\frac{1}{\lambda}\ln\frac{q_+}{x_0+\eps U(\infty)}+\frac{1}{\lambda}\ln \frac{x_0+\eps U(\infty)}{x_0+\eps U(\theta_\eps)}\\
         &=\frac{1}{\lambda}\ln q_++\frac{1}{\lambda}\ln\frac{1}{\eps}-\frac{1}{\lambda}\ln\left(U(\infty)+\frac{x_0}{\eps}\right)+\frac{1}{\lambda}\ln \frac{x_0+\eps U(\infty)}{x_0+\eps U(\theta_\eps)},
\end{align*}
so, introducing a standard Gaussian random variable $N=U(\infty)\sqrt{2\lambda}$, we can write
\begin{equation}
\label{eq:normalizing-theta}
\theta_\eps-\frac{2}{\lambda}\ln\frac{1}{\eps}=\frac{1}{\lambda}\ln (2\lambda q_+|x_0|)+\frac{1}{\lambda}R_\eps+Q_\eps,
\end{equation}
where
\[
R_\eps=-
\ln\left(N-\frac{|x_0|\sqrt{2\lambda}}{\eps}\right)
-\ln\left(\frac{|x_0|\sqrt{2\lambda}}{\eps}\right),
\]
and
\[
Q_\eps=\frac{1}{\lambda}\ln \frac{x_0+\eps U(\infty)}{x_0+\eps U(\theta_\eps)}.
\]
Let us recall that for each $\eps>0$, we are studying the above random variables on the event $F_\eps$ that
can be expressed as $F_\eps=\{N>|x_0|\sqrt{2\lambda}/\eps\}$.

\begin{lemma}\label{lem:log-conditioning} As $\eps\to0$,\quad $\Law[R_\eps|F_\eps]\Rightarrow \Lambda$.
\end{lemma}
\bpf There are several ways to prove this statement. 
One way is to introduce  $r(\eps)=|x_0|\sqrt{2\lambda}/\eps$ and directly write,
as in the proof of Corollary~1 of~\cite{Bakhtin-Gumbel},
\begin{align*}\notag
\frac{\Pp\{-\ln(N-r)-\ln r<x\}}{\Pp\{N>r\}}&=\frac{\Pp\{N>r+e^{-x}/r\}}{\Pp\{N>r\}}\\
\label{eq:proof_of_conv_to_gumbel}
&\sim\frac{
\frac{1}{\sqrt{2\pi}(r+e^{-x}/r)}e^{-(r+e^{-x}/r)^2/2}
}{
\frac{1}{\sqrt{2\pi}r}e^{-r^2/2}}
\sim \Lambda(x),\quad r\to\infty.
\notag
\end{align*}
\epf

\begin{remark}\rm
The statement of Lemma~\ref{lem:log-conditioning}  is not specific to the Gaussian distribution of the random variable $N$ as the proof may suggest.  Theorem~5 of~\cite{Bakhtin-Gumbel} shows that convergence of distributions of 
logarithms of residual life times to the Gumbel distribution holds if $N$ is replaced by any other distribution belonging to the domain
of attraction of the Gumbel distribution as a max-stable law. In fact, the theory of residual life times allows to describe the entire small
pool of distributions that can appear as a result of such a limiting procedure. We refer to~\cite{Bakhtin-Gumbel} for further explanation of connection
between extreme values, residual life times, and conditional exit times. We also notice that the proof of Theorem~\ref{th:Day}
in Section~\ref{sec:exit-to-0} is based on the same kind of asymptotics.
\end{remark}

\begin{lemma}\label{lem:Q_eps} As $\eps\to 0$,\quad $\Law[Q_\eps|F_\eps]\Rightarrow \delta_0$  (Dirac measure at $0$).
\end{lemma}
\bpf
We need to prove that $(x_0+\eps U(\infty))/(x_0+\eps U(\theta_\eps))$ given $F_\eps$, converges in probability to $1$, i.e.,
$\eps (U(\infty)-U(\theta_\eps))/(x_0+\eps U(\theta_\eps))$ converges to~$0$. Due to \eqref{eq:equation-for-theta},
it is sufficient to prove that, given $F_\eps$,
\[
\Delta_\eps=e^{\lambda\theta_\eps}\eps (U(\infty)-U(\theta_\eps))=e^{\lambda\theta_\eps}\eps\int_{\theta_\eps}^{\infty}e^{-\lambda s}dW(s)
\] 
converges to $0$. For that it is sufficient to check the
following $L^2$-convergence:
\[
\lim_{\eps\to 0}\frac{\E[ \Delta_\eps^2\ONE_{F_\eps}]}{\Pp(F_\eps)}=0,
\]
which, in turn, follows from
\begin{equation}\label{eq:L2-convergence}
\lim_{\eps\to 0}\frac{\E [\Delta_\eps^2\ONE_{D_\eps}]}{\Pp(F_\eps)}=0,
\end{equation}
where $D_\eps$ was introduced in~\eqref{eq:D_eps}, so $F_\eps\subset D_\eps$ since $x_0+\eps U(\infty)>0$
implies $x_0+\eps U(t)=0$ for some $t$. 

Let us now introduce $(\Fc_t)_{t\ge0}$, the completion of the natural filtration of the driving Wiener process $W$ (we implicitly have used it
above when dealing with the strong Markov property). 
Since $D_\eps\in\Fc_{\tau_\eps}$, we can use the It\^o isometry to derive
\begin{align*}
\E [\Delta_\eps^2\ONE_{D_\eps}]&=\eps^2 \E \left[\ONE_{D_\eps} e^{2\lambda\theta_\eps}\E\left[\left(\int_{\theta_\eps}^{\infty}e^{-\lambda s}dW(s)\right)^2\Bigr|\ \Fc_{\theta_\eps}\right]\right]\\
&=\eps^2 \E \left[\ONE_{D_\eps} e^{2\lambda\theta_\eps}\frac{e^{-2\lambda\theta_\eps}}{2\lambda}\right]=\frac{\eps^2}{2\lambda}\Pp(D_\eps),
\end{align*}
so~\eqref{eq:L2-convergence} follows since the reflection principle gives $\Pp(D_\eps)=2\Pp(F_\eps)$. \epf

Now, combining~\eqref{eq:normalizing-theta} with Lemmas~\ref{lem:log-conditioning} and~\ref{lem:Q_eps}, we obtain that 
\begin{equation}
\label{eq:convergence_under_F}
\Law\left[\theta_\eps-\frac{2}{\lambda}\ln\frac{1}{\eps}\Bigr|\ F_\eps\right]\Rightarrow\Law\left[\frac{\ln (2\lambda q_+|x_0|)}{\lambda}+\frac{1}{\lambda}Z\right],\quad \eps\to 0.
\end{equation}
Lemmas~\ref{lem:C-similar-to-E} and~\ref{lem:E-similar-to-F} allow to apply Lemma~\ref{lem:weak-convergence-under-equivalent-sequences}
and conclude that \eqref{eq:convergence_under_F} remains true if conditioning on~$F_\eps$ is replaced with
conditioning on~$C_\eps$. Recalling that $\theta_\eps=\tau_\eps$ on $C_\eps$, we complete the proof of~\eqref{eq:main-theorem-exit-right}.\epf

\section{Proof of Theorem~\ref{th:exit-scaling}}\label{sec:scaling}
 Since $X^1_\eps$
and $X^2_\eps$ in this model are disentangled, $\tau_\eps$ is essentially the exit time in the one-dimensional model studied above.

Using variation of constants, we can write
\[
X^2_\eps(\tau_\eps)=e^{-\mu \tau_\eps}X^2_\eps(0)+\eps e^{-\mu \tau_\eps}\int_0^{\tau_\eps} e^{\mu s}dW^2(s)=I_\eps+J_\eps.
\]
Let us study the limiting behavior of $I_\eps$ and $J_\eps$.  

Theorem~\ref{th:main} straightforwardly implies that
\begin{equation}
\label{eq:I_eps-convergence}
 \Law\left[I_\eps/\eps^{\frac{2\mu}{\lambda}+\alpha}|\ C_\eps\right]\Rightarrow \Law\left[(2\lambda q_+|x_1|)^{-\mu/\lambda} V^{\mu/\lambda }\xi\right],\quad\eps\to 0,
\end{equation}
where $V=e^{-Z}$ is a standard exponential random variable.

We can use the independence and Gaussianity of the processes involved  to
compute the characteristic function $\phi_\eps$ of $J_\eps/\eps$ conditioned on $C_\eps$:
\begin{align*}
\phi_\eps(r)&=\E\left[ e^{ir e^{-\mu \tau_\eps}\int_0^{\tau_\eps} e^{\mu s}dW^2(s)}\bigr| C_\eps\right]\\
&=
\E\left[ e^{-\frac{r^2e^{-2\mu \tau_\eps}}{2}\frac{e^{2\mu \tau_\eps}-1}{2\mu}}\Bigr| C_\eps\right]
=\E\left[ e^{-\frac{r^2}{2}\frac{1-e^{-2\mu \tau_\eps}}{2\mu}}\Bigr| C_\eps\right].
\end{align*} 
From Theorem~\ref{th:main} we know that,
conditioned on $C_\eps$,
$e^{-2\mu \tau_\eps}\to0$ in probability as $\eps\to 0$. Therefore, $\phi_\eps(r)\to e^{-r^2/(4\mu)}$ by dominated convergence,
i.e., 
\begin{equation}
\label{eq:J_eps-convergence}
\Law\bigl[J_\eps/\eps\ |\ C_\eps\bigr]\Rightarrow \Law\left[N/\sqrt{2\mu}\right],\quad\eps\to 0,
\end{equation}
 where $N$ is a standard Gaussian random variable.

In fact, joint convergence also holds true, and Theorem~\ref{th:exit-scaling} is an immediate corollary of the following result:
\begin{lemma}
\[
\Law\left[\left(\frac{I_\eps}{\eps^{\frac{2\mu}{\lambda}+\alpha}},\frac{J_\eps}{\eps}\right)\Bigr|\ C_\eps\right]\Rightarrow 
\Law\left[\left( \frac{V^{\mu/\lambda }\xi}{(2\lambda q_+|x_1|)^{\mu/\lambda}},\ \frac{N}{\sqrt{2\mu}}\right)\right],\quad\eps\to 0,
\]
where the joint distribution of $V,\xi,$ and $N$ is described in the statement of Theorem~\ref{th:exit-scaling}.
\end{lemma}
\bpf
We just have to slightly modify the proofs above. Let us compute $\phi_\eps(r_1,r_2)$, the characteristic function of 
$(I_\eps/\eps^{2\mu/\lambda+\alpha},J_\eps/\eps)$ conditioned on~$C_\eps$. 
Denoting conditional expectation given $C_\eps$ by $\E_{C_\eps}$, and representing
$I_\eps/\eps^{2\mu/\lambda+\alpha}=F_\eps(\xi_\eps,\tau_\eps)$ for some Borel function $F_\eps$, we can write
\begin{align*}
\phi_\eps(r_1,r_2)&=\E_{C_\eps} e^{i r_1F_\eps(\xi_\eps,\tau_\eps)}e^{ir_2 e^{-\mu \tau_\eps}\int_0^{\tau_\eps} e^{\mu s}dW^2(s)}
\\
&= \int_\R P_{\tau_\eps}(dt) \E_{C_\eps} e^{i r_1F_\eps(\xi_\eps,t)}e^{ir_2 e^{-\mu t}\int_0^{t} e^{\mu s}dW^2(s)},
\end{align*}
where $P_{\tau_\eps}(dt)$ denotes the distribution of $\tau_\eps$. We can use independence of $W^2$ from $C_\eps$ and $\xi_\eps$, and 
continue the computation:
\begin{align*}
\phi_\eps(r_1,r_2)&=\int_\R P_{\tau_\eps}(dt)\E_{C_\eps} e^{i r_1F_\eps(\xi_\eps,t)}\E_{C_\eps}e^{ir_2 e^{-\mu t}\int_0^{t} e^{\mu s}dW^2(s)}
               \\&=\int_\R P_{\tau_\eps}(dt)\E_{C_\eps} e^{i r_1F_\eps(\xi_\eps,t)}e^{-r_2^2\frac{1-e^{-2\mu t}}{4\mu}}
               \\&=\E_{C_\eps} e^{i r_1F_\eps(\xi_\eps,\tau_\eps)}e^{-r_2^2\frac{1-e^{-2\mu \tau_\eps}}{4\mu}}.
               \\&=\E_{C_\eps} e^{i r_1F_\eps(\xi_\eps,\tau_\eps)}e^{-r_2^2/(4\mu)}+\E_{C_\eps} e^{i r_1F_\eps(\xi_\eps,\tau_\eps)}
               (e^{r_2^2e^{-2\mu \tau_\eps}/(4\mu)}-1).
\end{align*}
The second term converges to 0 as the integral of a bounded random variable converging to zero in probability.
Convergence of the first term follows from~\eqref{eq:I_eps-convergence}.
\epf

\bibliographystyle{alpha}
\bibliography{gumbel}

\end{document}